\documentstyle[12pt]{article}
\catcode`@=11

\newskip\plaincentering \plaincentering=0pt plus 1000pt minus 1000pt
\def\@plainlign{\tabskip=0pt\everycr={}}
\def\eqalignno#1{\displ@y \tabskip\plaincentering
  \halign to\displaywidth{\hfil$\@lign\displaystyle{##}$\tabskip\z@skip
    &$\@lign\displaystyle{{}##}$\hfil\tabskip\plaincentering
    &\llap{$\@lign##$}\tabskip\z@skip\crcr
    #1\crcr}}
\def\leqalignno#1{\displ@y \tabskip\plaincentering
  \halign to\displaywidth{\hfil$\@lign\displaystyle{##}$\tabskip\z@skip
    &$\@lign\displaystyle{{}##}$\hfil\tabskip\plaincentering
    &\kern-\displaywidth\rlap{$\@lign##$}\tabskip\displaywidth\crcr
    #1\crcr}}
\def\plainLet@{\relax\iffalse{\fi\let\\=\cr\iffalse}\fi}
\def\plainvspace@{\def\vspace##1{\noalign{\vskip##1}}}

\def\intic@{\mathchoice{\hskip5\p@}{\hskip4\p@}{\hskip4\p@}{\hskip4\p@}}
\def\negintic@
 {\mathchoice{\hskip-5\p@}{\hskip-4\p@}{\hskip-4\p@}{\hskip-4\p@}}
\def\intkern@{\mathchoice{\!\!\!}{\!\!}{\!\!}{\!\!}}
\def\intdots@{\mathchoice{\cdots}{{\cdotp}\mkern1.5mu
    {\cdotp}\mkern1.5mu{\cdotp}}{{\cdotp}\mkern1mu{\cdotp}\mkern1mu
      {\cdotp}}{{\cdotp}\mkern1mu{\cdotp}\mkern1mu{\cdotp}}}
\newcount\intno@
\def\iint{\intno@=\tw@\futurelet\next\ints@}
\def\iiint{\intno@=\thr@@\futurelet\next\ints@}
\def\iiiint{\intno@=4 \futurelet\next\ints@}
\def\idotsint{\intno@=\z@\futurelet\next\ints@}
\def\ints@{\findlimits@\ints@@}
\newif\iflimtoken@
\newif\iflimits@
\def\findlimits@{\limtoken@false\limits@false\ifx\next\limits
 \limtoken@true\limits@true\else\ifx\next\nolimits\limtoken@true\limits@false
    \fi\fi}
\def\multintlimits@{\intop\ifnum\intno@=\z@\intdots@
  \else\intkern@\fi
    \ifnum\intno@>\tw@\intop\intkern@\fi
     \ifnum\intno@>\thr@@\intop\intkern@\fi\intop}
\def\multint@{\int\ifnum\intno@=\z@\intdots@\else\intkern@\fi
   \ifnum\intno@>\tw@\int\intkern@\fi
    \ifnum\intno@>\thr@@\int\intkern@\fi\int}
\def\ints@@{\iflimtoken@\def\ints@@@{\iflimits@
   \negintic@\mathop{\intic@\multintlimits@}\limits\else
    \multint@\nolimits\fi\eat@}\else
     \def\ints@@@{\multint@\nolimits}\fi\ints@@@}
\def\Sb{_\bgroup\vspace@
        \baselineskip=\fontdimen10 \scriptfont\tw@
        \advance\baselineskip by \fontdimen12 \scriptfont\tw@
        \lineskip=\thr@@\fontdimen8 \scriptfont\thr@@
        \lineskiplimit=\thr@@\fontdimen8 \scriptfont\thr@@
        \Let@\vbox\bgroup\halign\bgroup \hfil$\scriptstyle
            {##}$\hfil\cr}
\def\endSb{\crcr\egroup\egroup\egroup}
\def\Sp{^\bgroup\vspace@
        \baselineskip=\fontdimen10 \scriptfont\tw@
        \advance\baselineskip by \fontdimen12 \scriptfont\tw@
        \lineskip=\thr@@\fontdimen8 \scriptfont\thr@@
        \lineskiplimit=\thr@@\fontdimen8 \scriptfont\thr@@
        \Let@\vbox\bgroup\halign\bgroup \hfil$\scriptstyle
            {##}$\hfil\cr}
\def\endSp{\crcr\egroup\egroup\egroup}
\def\Let@{\relax\iffalse{\fi\let\\=\cr\iffalse}\fi}
\def\vspace@{\def\vspace##1{\noalign{\vskip##1 }}}
\def\aligned{\,\vcenter\bgroup\plainvspace@\plainLet@\openup\jot\m@th\ialign
  \bgroup \strut\hfil$\displaystyle{##}$&$\displaystyle{{}##}$\hfil\crcr}
\def\endaligned{\crcr\egroup\egroup}
\def\matrix{\,\vcenter\bgroup\plainLet@\plainvspace@
    \normalbaselines
  \m@th\ialign\bgroup\hfil$##$\hfil&&\quad\hfil$##$\hfil\crcr
    \mathstrut\crcr\noalign{\kern-\baselineskip}}
\def\endmatrix{\crcr\mathstrut\crcr\noalign{\kern-\baselineskip}\egroup
                \egroup\,}
\newtoks\hashtoks@
\hashtoks@={#}
\def\format{\crcr\egroup\iffalse{\fi\ifnum`}=0 \fi\format@}
\def\format@#1\\{\def\preamble@{#1}%
  \def\c{\hfil$\the\hashtoks@$\hfil}%
  \def\r{\hfil$\the\hashtoks@$}%
  \def\l{$\the\hashtoks@$\hfil}%
  \setbox\z@=\hbox{\xdef\Preamble@{\preamble@}}\ifnum`{=0 \fi\iffalse}\fi
   \ialign\bgroup\span\Preamble@\crcr}

\def\cases{\left\{\,\vcenter\bgroup\plainvspace@
     \normalbaselines\openup\jot\m@th
      \plainLet@\ialign\bgroup$\displaystyle{##}$\hfil&\quad$\displaystyle{{}##}$\hfil\crcr
      \mathstrut\crcr\noalign{\kern-\baselineskip}}
\def\endcases{\endmatrix\right.}
\newif\iftagsleft@
\tagsleft@true
\def\TagsOnRight{\global\tagsleft@false}
\def\tag#1$${\iftagsleft@\leqno\else\eqno\fi
 \hbox{\def\pagebreak{\global\postdisplaypenalty-\@M}%
 \def\nopagebreak{\global\postdisplaypenalty\@M}\rm(#1\unskip)}%
  $$\postdisplaypenalty\z@\ignorespaces}
\interdisplaylinepenalty=\@M
\def\plainallowdisplaybreak@{\def\allowdisplaybreak{\noalign{\allowbreak}}}
\def\plaindisplaybreak@{\def\displaybreak{\noalign{\break}}}
\def\align#1\endalign{\def\tag{&}\plainvspace@\plainallowdisplaybreak@\plaindisplaybreak@
  \iftagsleft@\plainlalign@#1\endalign\else
   \plainralign@#1\endalign\fi}
\def\plainralign@#1\endalign{\displ@y\plainLet@\tabskip\plaincentering\halign to\displaywidth
     {\hfil$\displaystyle{##}$\tabskip=\z@&$\displaystyle{{}##}$\hfil
       \tabskip=\plaincentering&\llap{\hbox{\rm(##\unskip)}}\tabskip\z@\crcr
             #1\crcr}}
\def\plainlalign@
 #1\endalign{\displ@y\plainLet@\tabskip\plaincentering\halign to \displaywidth
   {\hfil$\displaystyle{##}$\tabskip=\z@&$\displaystyle{{}##}$\hfil
   \tabskip=\plaincentering&\kern-\displaywidth
        \rlap{\hbox{\rm(##\unskip)}}\tabskip=\displaywidth\crcr
               #1\crcr}}
\def\proclaim#1{\par\bf #1. \unskip \it }
\def\endproclaim{\par\rm }
\def\demo#1{\par\it #1. \unskip \rm }
\def\enddemo{\par\rm }

\def\re@#1{\par\hangindent\parindent\indent\llap{#1\enspace}\ignorespaces}
\def\qfootnote#1{\edef\@sf{\spacefactor\the\spacefactor}{}#1\@sf
      \insert\footins{\let\egroup=}\footnotesize 
      \interlinepenalty100 \let\par=\endgraf
        \leftskip=0pt \rightskip=0pt
        \splittopskip=10pt plus 1pt minus 1pt \floatingpenalty=20000
   \smallskip\re@{#1}\bgroup\strut\aftergroup{\strut\egroup}\let\next}
\catcode`\@=\active
\TagsOnRight

\begin{document}
\begin{center}
{\large\bf
 On $L^{1}$ Convergence of Fourier Series
 \\ of Complex Valued
 Functions}

\vspace{5mm}
 R. J. Le$^{1,2}$ and S. P. Zhou$^{1}$
\footnote{ 1. Institute of Mathematics, Zhejiang Sci-Tech
University, Xiasha Ecnomic Development Area, Hangzhou, Zhejiang
310018 China. 2. Department of Mathematics, Ningbo University, Ningbo, Zhejiang 315211  China
  Supported in part by    Natural Science Foundation of China under grant number 10471130. \\

  The main part of the work contained in this paper was done while the first named author was working with the second named author at Institute of Mathematics, Zhejiang Sci-Tech
  University.\\

Key wards and phrases: $L^{1}$-convergence, Fourier series,
monotonicity, convexity, complex valued function}

\end{center}

 \abstract In the present
paper, we give a brief review of $L^{1}$-convergence of
trigonometric series. Previous known results in this direction are
improved and generalized by establishing a new condition.
 \endabstract

\small 1991 Mathematics Subject Classification. 42A20 42A32\\

\bf \S 1. Introduction\newline\newline

\rm Let $L_{2\pi}$ be the space of all real or complex valued integrable functions $f(x)$ of period $2\pi$ with norm
$$\|f\| = \int_{-\pi}^{\pi}|f(x)|dx.$$
When $f(x)\in L_{2\pi}$ is a real valued even function, we denote the Fourier series of $f$ by
$$\frac{a_{0}}{2}+\sum_{k=1}^{\infty}a_{k}\cos kx,\tag{1}$$
and its $n$th partial sum $S_{n}(f,x)$ by
$$\frac{a_{0}}{2}+\sum_{k=1}^{n}a_{k}\cos kx.$$

There is a long history for $L^{1}$-convergence of Fourier series or trigonometric series. The story   started from monotonicity of the coefficients, the original theorem can be stated as follows:
\proclaim{Theorem A}
 Let $f(x)\in L_{2\pi}$ with Fourier series $(1)$.  If $\{a_{n}\}$ is a nonincreasing sequence, then
$$\lim_{n\to\infty}\|f-S_{n}(f)\|=0\;\; {\rm if and only if}\;\;\lim_{n\to\infty}a_{n}\log n=0.$$
\endproclaim

The course of development for this way can be summed as follows:

\bf Condition: $f\in L_{2\pi}$

$${\rm (coefficients)}\;\;{\rm nonincreasing}\;\;\Rightarrow\;\;{\rm quasimonotone}$$
$$\Rightarrow\;\;{\rm regularly varying quasimonotone}$$
$$\Rightarrow\;\;O-{\rm regularly varying quasimonotone}$$

Conclusion: \rm $\lim\limits_{n\to\infty}\|f-S_{n}(f)\|=0$  if and only if $\lim\limits_{n\to\infty}a_{n}\log n=0.$

Interested readers could check references such as [1, 4-10].

For a sequence $\{c_{n}\}_{n=0}^{\infty}$, let
$$\Delta c_{n} = c_{n}-c_{n+1}.$$

A non-decreasing positive sequence $\{R(n)\}_{n=1}^{\infty}$ is said to be $O$-regularly varying if\footnote{In some papers, this requirement is written as that for some $\lambda>1$, $\limsup\limits_{n\to\infty}R([\lambda n])/R(n)<\infty$. One can easily check out that they are equivalent.}
$$\limsup_{n\to\infty}\frac{R(2n)}{R(n)}<\infty.$$

A complex sequence $\{c_{n}\}_{n=0}^{\infty}$, with ${\rm
Re}c_{n}\geq 0$, is $O$-regularly varying quasimonotone in complex
sense if for some $\theta_{0}\in [0,\pi/2)$ and some $O$-regularly
varying sequence $\{R(n)\}$ the sequence
$$\Delta\frac{c_{n}}{R(n)}\in K(\theta_{0}):=\{z: |{\rm arg} z|\leq\theta_{0}\},\;\;n=1, 2,\cdots.$$

Evidently, if $\{c_{n}\}$ is a real sequence, then the $O$-regularly varying quasimonotonicity becomes
$$\Delta\frac{c_{n}}{R(n)}\geq 0,\;\;n=1, 2, \cdots,$$
which was used in many works to generalize the regularly varying quasimonotone condition and, in particular, the quasimonotone condition\footnote{A real sequence $\{b_{n}\}_{n=0}^{\infty}$
is defined to be quasimonotone if, for some $\alpha\geq 0$, the sequence $\{b_{n}/n^{\alpha}\}$ is non-increasing.}.

We give a new condition to unify and generalize various quasimonotonicity conditions.

\proclaim{Definition}
 Let {\bf c}$=\{c_{n}\}_{n=1}^{\infty}$ be a sequence satisfying $c_{n}\in K(\theta_{1})$ for some $\theta_{1}\in [0, \pi/2)$ and $n=1, 2, \cdots$. If there is a natural number $N_{0}$ such that
$$\sum_{n=m}^{2m}|\Delta c_{n}|\leq M({\bf c})\max_{m\leq n<m+N_{0}}|c_{n}|\tag2$$
holds for all $m=1, 2, \cdots$, where  $M({\bf c})$ indicates a
positive constant only depending upon ${\bf c}$, then we say that
the sequence ${\bf c}$ belongs to class {\rm GBV}.
\endproclaim

\rm We recall the following results.

\proclaim{Lemma XZ}
 Suppose a complex sequence $\{c_{n}\}$ is $O$-regularly varying quasimonotone,
then there is a positive constant $M$ depending upon $\theta_{0}$ only such that
$$|c_{n}| \leq M{\rm Re}c_{n},\;\;n=1, 2,\cdots,$$
or in other words, $c_{n}\in K(\theta_{1})$ for some $\theta_{1}\in [0, \pi/2)$ and $n=1, 2, \cdots$.
\endproclaim

   \rm The argument exactly follows from Xie and Zhou [10, Lemma 1], and the condition $\lim\limits_{n\to\infty}c_{n}=0$ there can be cancelled however.

\proclaim{Lemma ZL} Let $\{c_{n}\}$ be any given complex $O$-regularly varying quasimonotone sequence. Then $\{c_{n}\}$ satisfies $(2)$ for $N_{0}=1$.
\endproclaim

 \rm The proof can be copied from Le and Zhou [11, Theorem 3] with omitting the condition $\lim\limits_{n\to\infty}c_{n}=0$ there.

We clearly see that, from Lemma XZ and Lemma ZL, if ${\bf
c}=\{c_{n}\}_{n=1}^{\infty}$ is any given complex $O$-regularly
varying quasimonotone sequence, then $\{c_{n}\}\in {\rm GBV}$. The
converse is obviously not true. Thus the class {\rm GBV} gives an
essential and explicit generalization to the class of
$O$-regularly varying quasimonotone sequences.

On the other hand, Leindler in [2] introduced the following ``rest bounded variation" sequences (RBVS sequences): Let \mbox{\bf b}$=\{b_{n}\}_{n=1}^{\infty}$ be a nonnegative sequence satisfying $\lim\limits_{n\to\infty}b_{n}=0$ and
$$\sum_{n=m}^{\infty}|b_{n}-b_{n+1}|\leq M(\mbox{\bf b})b_{m}$$
for some constant $M(\mbox{\bf b})$ depending only upon $\mbox{\bf
b}$ and $m=1, 2, \cdots$. Also, Leindler [3, Theorem 1] indicated
that quasimonotonicity  and the above ``rest bounded variation"
condition are not comparable. The idea of {\rm GBV} sequences
gives a tricky generalization of RBVS sequences.

Let $f(x)\in L_{2\pi}$ be a complex valued function,
denote the Fourier series of $f$ by
$$\sum_{k=-\infty}^{\infty}\hat{f}(k)e^{ikx},$$
and its $n$th partial sum by
$$\sum_{k=-n}^{n}\hat{f}(k)e^{ikx}.$$

In Section 2, we will establish Theorem 1, which contains and generalizes all results we stated in the beginning.

\proclaim{Theorem 1}
 Let $f(x)\in L_{2\pi}$ be a complex valued function. If the Fourier coefficients of $f$ satisfy that
$$\lim_{\lambda\to 1+0}\limsup_{n\to\infty}\sum_{k=n}^{[\lambda n]}|\Delta\hat{f}(k)-\Delta\hat{f}(-k)|\log k = 0,\tag3$$
and $\{\hat{f}(n)\}_{n=0}^{+\infty}\in {\rm GBV}$, then\\
$$\lim\limits_{n\to\infty}\|f-S_{n}(f)\|=0\;\;\mbox{\rm  if  and  only
if}\;\;\lim\limits_{n\to\infty}\hat{f}(n)\log|n|=0.$$
\endproclaim

\rm We make a quick remark here. We notice that, in the complex
valued function space, the assumption that both
$\{\hat{f}(n)\}_{n=0}^{\infty}\in {\rm GBV}$ and
$\{\hat{f}(-n)\}_{n=1}^{\infty}\in {\rm GBV}$   is a convenient
one, but is almost trivial, since it is almost the same as the
condition in real case. Thus people usually use one side condition
with some kind of balance conditions in considering those
problems. Theorem 1 reflects this kind of thinking.

\rm From $\sum_{k=1}^{\infty}k^{-\alpha}\sin 2^{k^{2}}x=:\sum_{n=1}^{\infty}b_{n}\sin nx$, $1<\alpha<2$, we can clearly see that,  $b_{n}\log n\to \infty$, $n\to\infty$, therefore, the condition $(2)$, in general sense, cannot be further generalized.

On the basis of the proof of Theorem 1, we will establish some corresponding $L^{1}$-approximation theorems, particularly the following

\proclaim{Corollary 3}
 Let $f(x)\in L_{2\pi}$ be a complex valued function,
and both $\{\hat{f}(n)\}_{n=0}^{+\infty}\in {\rm GBV}$ and
$\{\hat{f}(-n)\}_{n=1}^{+\infty}\in {\rm GBV}$. If $f(x)$ has
$r-1$ absolutely continuous derivatives, then
$$\|f-S_{n}(f)\|=O\left((n+1)^{-r}\omega\left(f^{(r)}, (n+1)^{-1}\right)_{L}\right)$$
if and only if
$$\hat{f}(n)\log|n|
=O\left((|n|+1)^{-r}\omega\left(f^{(r)}, (|n|+1)^{-1}\right)_{L}\right),$$
where $\omega(f,t)_{L}$ is the modulus of continuity of $f$ in integral metric, that is,
$$\omega(f,t)_{L} = \max_{0\leq h\leq t}\|f(x+h)-f(x)\|.$$
\endproclaim

Throughout the paper, $C$ denotes a positive constant (which is
independent of $n$ and $x\in [0,2\pi]$) not necessarily the same
at each occurrence. In some specific cases, we also use $M({\bf
c})$ to indicate a positive constant only depending upon the
sequence ${\bf c}$.\newline\newline

 \bf\S 2. Convergence and Approximation \newline\newline

\rm We need to establish several lemmas first.

\proclaim{Lemma 1} Let $\{c_{n}\}\in{\rm GBV}$, then there is a
positive constant $M$ depending upon ${\bf c}$ and $\theta_{1}$
only such that, for $j=0, 1, \cdots, [n/N_{0}]-1$,
$$|c_{2n}|\leq M({\bf c}, \theta_{1})\left(\max_{n+jN_{0}\leq k< n+(j+1)N_{0}}{\rm Re}c_{k}+{\rm Re}c_{2n+2jN_{0}}\right).$$
\endproclaim

\demo{Proof}  For $j=0, 1, \cdots, [n/N_{0}]-1$, we see that
$$n\leq n+jN_{0}\leq n+(j+1)N_{0}\leq 2n.$$
From the conditions of GBV, we get for $0\leq j\leq [n/N_{0}]-1$,
$$|c_{2n}|=\left|\sum_{k=2n}^{2n+2jN_{0}-1}\Delta c_{k}+ c_{2n+2jN_{0}}\right|\leq \sum_{k=n+jN_{0}}^{2n+2jN_{0}}|\Delta c_{k}|+|c_{2n+2jN_{0}}|$$
$$\leq M({\bf c})\max_{n+jN_{0}\leq k< n+(j+1)N_{0}}|c_{k}|+|c_{2n+2jN_{0}}|$$
$$\leq M({\bf c}, \theta_{1})\left(\max_{n+jN_{0}\leq k<n+(j+1)N_{0}}{\rm Re}c_{k}+{\rm Re}c_{2n+2jN_{0}}\right).
$$
\enddemo

\proclaim{Lemma 2} Let $\{c_{n}\}\in{\rm GBV}$, then for given
$\lambda>1$, one has
$$\sum_{k=n}^{[\lambda n]}|\Delta c_{k}|\log k = O\left(\max_{n\leq k\leq [\lambda n]}|c_{k}|\log k\right), \;\;n=1,2,\cdots.$$
\endproclaim

\demo{Proof} For large enough $n$, from condition (2), it is
simply a straightforward result.
\enddemo

\proclaim{Lemma 3} Write

$$\phi_{\pm n}(x) = \sum_{k=1}^{n}\frac{1}{k}\left(e^{i(k\mp n)x}-e^{-i(k\pm n)x}\right),$$
then

$$|\phi_{\pm n}(x)| \leq 6\sqrt{\pi}.$$
\endproclaim

\demo{Proof} It is a direct consequence of the well-known inequality

$$\sup_{n\geq 1}\left|\sum_{k=1}^{n}\frac{\sin kx}{k}\right|\leq 3\sqrt{\pi}.$$
\enddemo

\demo{Proof of Theorem 1} Given $\epsilon>0$, by (3), there is a $\lambda>1$ such that
$$\sum_{k=n}^{[\lambda n]}\left|\Delta\hat{f}(k)-\Delta\hat{f}(-k)\right|\log k \leq \epsilon\tag5$$
holds for all sufficiently large $n>0$.  Denote
$$\tau_{\lambda n,n}(f,x) = \frac{1}{[\lambda n]-n}\sum_{k=n}^{[\lambda n]-1}S_{k}(f,x),$$
then clearly,
$$\lim_{n\to\infty}\|f-\tau_{\lambda n,n}(f)\| = 0.\tag6$$
Write
$$D_{k}(x) = \frac{\sin((2k+1)x/2)}{2\sin(x/2)},$$
$$D_{k}^{\ast}(x) = \cases
\frac{\cos(x/2)-\cos((2k+1)x/2)}{2\sin(x/2)},& |x|\leq 1/n,\\
-\frac{\cos((2k+1)x/2)}{2\sin(x/2)},&1/n\leq |x|\leq \pi,\endcases$$
$$E_{k}(x) = D_{k}(x) + iD_{k}^{\ast}(x).$$
For $k=n, n+1, \cdots, 2n$, we have
$$E_{k}(\pm x)-E_{k-1}(\pm x) = e^{\pm ikx},\tag7$$
$$E_{k}(x)+E_{k}(-x) = 2D_{k}(x),\tag8$$
and
$$\|E_{k}\|+\|D_{k}\| = O(\log k).\tag9$$
By using (7), (8) and applying Abel Transformation we get
$$\tau_{\lambda n,n}(f,x)-S_{n}(f,x) = \frac{1}{[\lambda n]-n }\sum_{k=n+1}^{[\lambda n]}([\lambda n]-k)\left(\hat{f}(k)e^{ikx}+\hat{f}(-k)e^{-ikx}\right)\tag10$$
$$= \frac{1}{[\lambda n]-n }\sum_{k=n}^{[\lambda n]}([\lambda n]-k)\left(2\Delta\hat{f}(k)D_{k}(x)-(\Delta\hat{f}(k)-\Delta\hat{f}(-k))E_{k}(-x)\right)$$
$$+ \frac{1}{[\lambda n]-n }\sum_{k=n}^{[\lambda n]-1}\left(\hat{f}(k+1)E_{k}(x)-\hat{f}(-k-1))E_{k}(-x)\right)$$
$$- \left(\hat{f}(n)E_{n}(x)+\hat{f}(-n)E_{n}(-x)\right).$$
Thus (9) yields that
$$\|\tau_{\lambda n,n}(f)-S_{n}(f)\| = O\left(\sum_{k=n}^{[\lambda n]}\left|\Delta\hat{f}(k)\right|\log k\right) + O\left(\sum_{k=n}^{[\lambda n]}\left|\Delta\hat{f}(k)-\Delta\hat{f}(-k)\right|\log k\right)$$
$$+ O\left(\max_{n\leq|k|\leq [\lambda n]-1}\left|\hat{f}(k)\right|\log|k|\right).$$
With Lemma 2, one has
$$\|f-S_{n}(f)\| \leq \|f-\tau_{\lambda n,n}(f)\|+O\left( \sum_{k=n}^{[\lambda n]}\left|\Delta\hat{f}(k)-\Delta\hat{f}(-k)\right|\log k\right)$$
$$+ O\left(\max_{n\leq|k|\leq [\lambda n]}\left|\hat{f}(k)\right|\log|k|\right),$$
then
$$\limsup_{n\to\infty}\|f-S_{n}(f)\|\leq\epsilon$$
follows from (5), (6) and
$$\lim_{n\to\infty}\hat{f}(n)\log|n|=0,\tag11$$
and we have obtained that
$$\lim_{n\to\infty}\|f-S_{n}(f)\|=0.\tag12$$

Now we come to prove that (12) implies (11). From Lemma 3, we deduce that
$$\frac{1}{6\sqrt{\pi}}\left|\int_{-\pi}^{\pi}(f(x)-S_{n}(f, x))\phi_{n}(x)dx\right|\leq\|f(x)-S_{n}(f, x)\|,$$
so that
$$\sum_{k=1}^{n}\frac{1}{k}\hat{f}(n+k)=O(\|f(x)-S_{n}(f, x)\|),$$
in particular,
$$\sum_{k=1}^{n}\frac{1}{k}{\rm Re}\hat{f}(n+k)=O(\|f(x)-S_{n}(f, x)\|).$$
By the same technique, we also have
$$\sum_{k=1}^{2n}\frac{1}{k}{\rm Re}\hat{f}(2n+k)=O(\|f(x)-S_{2n}(f, x)\|).$$
Applying Lemma 1, we obtain that
$$|\hat{f}(2n)|\log n\leq M(N_{0})|\hat{f}(2n)|\sum_{j=1}^{[n/N_{0}]-1}j^{-1}$$
$$\leq M({\bf c}, \theta_{1}, N_{0})\left(\sum_{j=1}^{[n/N_{0}]-1}(N_{0}j)^{-1}\left(\max_{jN_{0}\leq k<(j+1)N_{0}}{\rm Re}\hat{f}(n+k)\right.\right.$$
$$\left.+ \sum_{j=1}^{[n/N_{0}]-1}(N_{0}j)^{-1}{\rm Re}\hat{f}(2n+2jN_{0})\right)$$
 $$\leq M({\bf c}, \theta_{1}, N_{0})\left(\sum_{k=1}^{n}\frac{1}{k}{\rm Re}\hat{f}(n+k)+ \sum_{k=1}^{2n}\frac{1}{k}{\rm Re}\hat{f}(2n+k)\right)$$
$$\leq M({\bf c}, \theta_{1}, N_{0})\left(\|f(x)-S_{n}(f, x)\|+ \|f(x)-S_{2n}(f, x)\|\right),$$
accordingly (6) implies that
$$\lim_{n\to +\infty}|\hat{f}(2n)|\log n = 0.\tag13$$
The same argument can be applied as well to achieve that
$$\lim_{n\to +\infty}|\hat{f}(2n+1)|\log n = 0.\tag14$$
Recalling (10) and (9), we see that for the given $\lambda>1$,
$$\|\hat{f}(-n)E_{n}(-x)\|\leq \|\tau_{\lambda n,n}(f)-S_{n}(f)\| + \frac{1}{[\lambda n]-n }\left\|\sum_{k=n}^{[\lambda n]-1}\hat{f}(-k-1)E_{k}(-x)\right\| $$
$$+ O\left(\sum_{k=n}^{[\lambda n]}\left(|\Delta\hat{f}(k)-\Delta\hat{f}(-k)|\log k+
\left|\Delta\hat{f}(k)\right|\log k\right) + \max_{n\leq k\leq [\lambda n]}|\hat{f}(k)|\log k\right). \tag15 $$
Set
$$I :=
\int_{n^{-1}\leq |x|\leq\pi}\left|\frac{1}{2\sin(x/2)}\sum_{k=n}^{[\lambda n]-1}\hat{f}(-k-1)\left(\sin\frac{(2k+1)x}{2}+i\cos\frac{(2k+1)x}{2}\right)\right|dx.$$
It is easy to see that
$$\left\|\sum_{k=n}^{[\lambda n]-1}\hat{f}(-k-1)E_{k}(-x)\right\|
= I + O\left(n\max_{n<k\leq [\lambda n]}|\hat{f}(-k)|\right).$$
Since
$$I\leq\left(\int_{n^{-1}\leq|x|\leq\pi}\left|\sum_{k=n}^{[\lambda n]-1}\hat{f}(-k-1)\left(\sin\frac{(2k+1)x}{2}+i\cos\frac{(2k+1)x}{2}\right)\right|^{2}dx\right)^{1/2}$$
$$\times\sqrt{\int_{n^{-1}}^{\pi}\frac{2dx}{\sin^{2}(x/2)}},$$
and trigonometric function system is orthonormal, we thus have
$$I = O\left(\sqrt{n}\left(\sum_{k=n+1}^{[\lambda n]}|\hat{f}(-k)|^{2}\right)^{1/2}\right) =O\left(n\max_{n<k\leq [\lambda n]}|\hat{f}(-k)|\right),$$
which implies that,
$$\frac{1}{[\lambda n]-n }\left\|\sum_{k=n}^{[\lambda n]-1}\hat{f}(-k-1)E_{k}(-x)\right\| = O\left(\max_{n<k\leq [\lambda n]}|\hat{f}(-k)|\right).\tag16$$
Combining (15) with (6), (12)-(14), (16) and Lemma 2 and noting that $\hat{f}(-n)\to 0$, $n\to\infty$ (since $f\in L_{2\pi}$), we have
$$\|\hat{f}(-n)E_{n}(-x)\|
\leq\sum_{k=n}^{[\lambda n]}|\Delta\hat{f}(k)-\Delta\hat{f}(-k)|\log k
 +O\left(\|\tau_{\lambda n, n}(f)-S_{n}(f)\|\right)\tag17$$
  $$+O\left(\max_{n\leq k\leq [\lambda n]}\left|\hat{f}(k)\right|\log k\right)
+O\left(\max_{n<k\leq [\lambda n]}|\hat{f}(-k)|\right)$$
$$\leq\sum_{k=n}^{[\lambda n]}|\Delta\hat{f}(k)-\Delta\hat{f}(-k)|\log k + o(1), \;\;n\to\infty.$$
At the same time, evidently we have
$$\|\hat{f}(-n)E_{n}(-x)\|\geq |\hat{f}(-n)|\|D_{n}(x)\|\geq \frac{1}{\pi}|\hat{f}(-n)|\log n.\tag18$$
Altogether (17), (18) and (5) yield that
$$|\hat{f}(-n)|\log n \leq \sum_{k=n}^{[\lambda n]}|\Delta\hat{f}(k)-\Delta\hat{f}(-k)|\log k\leq\epsilon$$
holds for sufficiently large $n$, which, together with (13) and
(14), completes the proof of (11). So far, we have proved Theorem
1.
\enddemo

According to Lemma 2, by the same argument, we see that condition (3) can be replaced by the following condition
$$\lim_{\lambda\to 1+0}\limsup_{n\to\infty}\sum_{k=n}^{[\lambda n]}|\Delta\hat{f}(-k)|\log k = 0,$$
while the same coclusion of Theorem 1 still holds. Therefore in a similar way we have
\proclaim{Corollary 1}
 Let $f(x)\in L_{2\pi}$ be a complex valued function.
If both $\{\hat{f}(n)\}_{n=0}^{+\infty}\in {\rm GBV}$ and
$\{\hat{f}(-n)\}_{n=0}^{+\infty}\in {\rm GBV}$, then
$$\lim_{n\to\infty}\|f-S_{n}(f)\|=0\;\; {\rm if and only if}\;\;\lim_{n\to\infty}\hat{f}(n)\log|n|=0.$$
\endproclaim

 If $f(x)$ is  a real valued function, then its Fourier coefficints $\hat{f}(n)$ and $\hat{f}(-n)$ are a pair of conjugate numbers. Consequently, $\{\hat{f}(n)\}_{n=0}^{+\infty}\in {\rm GBV}$ if and only if $\{\hat{f}(-n)\}_{n=0}^{+\infty}\in {\rm GBV}$. That is,

\proclaim{Corollary 2}
 Let $f(x)\in L_{2\pi}$ be a real valued function.
If $\{\hat{f}(n)\}_{n=0}^{+\infty}\in {\rm GBV}$, then
$$\lim_{n\to\infty}\|f-S_{n}(f)\|=0\;\; {\rm if and only if}\;\;\lim_{n\to\infty}\hat{f}(n)\log|n|=0.$$
\endproclaim

 \rm Let $E_{n}(f)_{L}$ be the best approximation of a complex valued function $f\in L_{2\pi}$ by trigonometric polynomials of degree $n$ in integral metric, that is,

$$E_{n}(f)_{L} = \inf_{c_{k}}\left\|f-\sum_{k=-n}^{n}c_{k}e^{ikx}\right\|.$$
We can, with omitting the proof, deduce the following $L^{1}$-approximation theorem in a similar way to Theorem 1.

\proclaim{Theorem 2}
 Let $f(x)\in L_{2\pi}$ be a complex valued function, $\{\psi_{n}\}$ is a decreasing sequence tending to zero with that
  $$\psi_{n}=O(\psi_{2n}).$$
If both $\{\hat{f}(n)\}_{n=0}^{+\infty}\in {\rm GBV}$ and
$\{\hat{f}(-n)\}_{n=1}^{+\infty}\in {\rm GBV}$, then
$$\|f-S_{n}(f)\|=O(\psi_{n})$$
if and only if
$$E_{n}(f)_{L}=O(\psi_{n})\;\;{\rm and}\;\;\hat{f}(n)\log|n|=O(\psi_{|n|}).$$
\endproclaim

Corollary 3 we stated in Section 1 is a direct consequence from Theorem 2.\newline\newline

\vspace{3mm}
\newpage
\begin{center}
{\Large\bf References}
\end{center}
\begin{enumerate}
\rm

\item W. O. Bary and \v{C}. V. Stanojevic, \it On weighted integrability of trigonometric series and $L^{1}$-convergence of Fourier series, \rm Proc. Amer. Math. Soc. 96(1986),53-61.
\item J. Karamata, \it Sur un mode de croissance r\'{e}guliere des fonctions, \rm Mathematica (Cluj) 4(1930),38-53.
\item  L. Leindler, \it On the uniform convergence and boundedness of a certain class of sine series, \rm Anal. Math. 27(2001), 279-285.
\item L. Leindler, \it A new class of numerical sequences and its applications to sine and cosine series, \rm Anal. Math. 28(2002), 279-286.
\item V. B. Stanojevic, \it $L^{1}$-convergence of Fourier series with complex quasimonotone coefficients, \rm Proc. Amer. Math. Soc. 86(1982), 241-247.
\item V. B. Stanojevic, \it Convergence of Fourier series with complex quasimonotone coefficients of bounded variation of order $m$, \rm J. Math. Anal. Appl. 115(1986), 482-505.
\item V. B. Stanojevic, \it $L^{1}$-convergence of Fourier series with complex quasimonotone coefficients, \rm Acad. Serbe Sci. Arts Glas 346(1986),29-48.
\item V. B. Stanojevic, \it $L^{1}$-convergence of Fourier series with $O$-regularly varying quasimonotone coefficients, \rm J. Approx. Theory 60(1990), 168-173.
\item S. A. Telyakovskii and G. A. Fomin, \it On the convergence in the $L$ metric of Fourier series with quasi-monotone coefficients, \rm Proc. Steklov Inst. Math. 134(1975), 351-355.
\item T. F. Xie and S. P. Zhou, \it $L^{1}$-approximation of Fourier series of complex valued functions, \rm Proc. Royal Soc. Edinburgh 126A(1996), 343-353.
\item  R. J. Le and S. P. Zhou,\it A new condition for the uniform convergence of certain trigonometric series, \rm Acta Math. Hungar. 108(2005), 161-169.
 \item A. Zygmund, \it Trigonometric Series, \rm Cambridge Univ. Press, Cambridge, 1959.
 \end{enumerate}
 \enddocument